# Double Reweighted Estimators for the Parameters of the Multivariate t Distribution

Fatma Zehra Doğru[1]*, Y. Murat Bulut[2] and Olcay Arslan[3]

[1] Giresun University, Department of Econometrics, 28100 Giresun/Turkey.
(E-mail: fatma.dogru@giresun.edu.tr)
[2] Eskisehir Osmangazi University, Department of Statistics, 26480 Eskisehir/Turkey.(E-mail: ymbulut@ogu.edu.tr)
[3] Ankara University, Department of Statistics, 06100 Ankara/Turkey.
(E-mail: oarslan@ankara.edu.tr)

**Abstract**

The t-distribution has many useful applications in robust statistical analysis. The parameter estimation of the t-distribution is carried out using ML estimation method, and the ML estimates are obtained via the EM algorithm. In this study, we consider an alternative estimation method for all the parameters of the multivariate-t distribution using the ML$q$ estimation method. We adapt the EM algorithm to obtain the ML$q$ estimates for all the parameters. We provide a small simulation study to illustrate the performance of the ML$q$ estimators over the ML estimators and observe that the ML$q$ estimators have considerable superiority over the ML estimators.

## 1. Introduction

The Shannon's entropy is defined as $H(X) = -E[\log p(X)]$, where $p(x)$ is the probability density function (pdf) of $X$. It was presented by Akaike (1973) that in a parametric model the maximization of the log-likelihood function is equivalent to the minimization of the empirical version of the Shannon's entropy $(-\sum_{i=1}^{n} \log p(X_i))$. Havrda and Charvát (1967) and Tsallis (1988) have introduced a $q$-extension of the Shannon entropy by using the $q$-logaritmic function which is given below

$$L_q(u) = \begin{cases} \log(u), & u \geq 0, q = 1, \\ \dfrac{u^{1-q} - 1}{1 - q}, & u \geq 0, q \neq 1. \end{cases} \quad (1)$$

Using the $q$-logaritmic function given in (1), Havrda- Charvát-Tsallis entropy, which is also called the $q$-entropy, is defined as

$$H_q(X) = -E[L_q(p(X))]. \quad (2)$$

Note that when $q = 1$ in (2), we get the Shannon entropy.

Using the empirical version of the $q$-entropy functional $H_q(X)$ in (2), Ferrari and Yang (2010) introduced the maximum L$q$-likelihood (ML$q$) estimation method and studied its properties. In the paper by Ferrari and Vecchia (2012), the authors considered on the robustness properties of the ML$q$ estimation and given the relationship between L$q$-likelihood estimation and the estimation by minimization of power divergences proposed by Basu et al. (1998). ML$q$ estimation method can yield robust estimators for moderate or small sample sizes which can provide an important progress with regards to mean squared error at the expense of a slight increase in bias. Note that the ML$q$ estimators can be considered as weighted likelihood estimators with the weights related to the $(1 - q)$ power to the density function.

It is well known that the t-distribution, which belongs to the family of elliptical distribution, is often used as a robust alternative distribution to the normal distribution. The maximum likelihood (ML) estimation method is one way to find the estimators for the parameters of the t distribution. The ML estimators are very popular in robust statistical analysis since they provide alternative robust estimators to the classical estimators obtained from the normal distribution. Since the likelihood equations cannot be solved analytically the Expectation-Maximization (EM) algorithm (Dempster et al. (1977), see also McLachlan and Krishnan (1997)) is one way to find the ML estimates for parameters of the t distribution.

However, it is well-known that when estimating the degrees of freedom parameter along with the other parameters the estimators become no longer locally robust due to the unboundedness of the score function for the degrees of freedom parameter. Therefore, to obtain robust estimators the degrees of freedom parameter is usually assumed to be known and considered as a robustness tuning constant (e.g., see, Lucas (1997) and Lange et al. (1989)). The purpose of this paper is to use the ML$q$ estimation as an alternative to the ML estimation method for the parameters of the multivariate *t* distribution. We observe that the score function for the degrees of freedom parameter obtained from the ML$q$ estimation method is bounded and so the proposed estimators for all the parameters will be locally robust unlike the estimators obtained from ML estimation method. We also noticed that the ML$q$ estimators, as we mentioned before, are the adaptively weighted form of the sample mean and the sample covariance matrix similar to the ML estimators. However, the weights in this case are faster decreasing than the weights obtained from the ML estimation method. This yields more robust estimators than the ML estimators in terms of down-weighting outlying observations. Concerning the computation of the ML$q$ estimates, we adapt the EM algorithm to obtain the ML$q$ estimates.

The rest of the paper is organized as follows. In Section 2, we briefly summarize some properties of multivariate t distribution. In Section 3, we give the ML estimation for multivariate t distribution. Also, we describe the ML$q$ estimation method and provide the ML$q$ estimators for the parameters of the multivariate t-distribution. In Section 4, we give some numerical examples to evaluate the performance of the ML$q$ estimators over the ML estimators. The paper is finalized with a conclusion section.

**2. Multivariate t-distribution**

Let $\boldsymbol{X} \in R^p$, be a $p$-dimensional random vector from multivariate t-distribution ($\boldsymbol{X} \sim t_p(\boldsymbol{\mu}, \Sigma, \nu)$). The probability density function (pdf) of $\boldsymbol{X}$ is given below

$$f(x; \boldsymbol{\mu}, \Sigma, \nu) = \frac{\Gamma\left(\frac{\nu + p}{2}\right) |\Sigma|^{-\frac{1}{2}}}{(\pi\nu)^{p/2} \Gamma\left(\frac{\nu}{2}\right)} \left[1 + \frac{1}{\nu}(x - \boldsymbol{\mu})^T \Sigma^{-1}(x - \boldsymbol{\mu})\right]^{-\frac{(\nu+p)}{2}}, \qquad (3)$$

where $\boldsymbol{\mu}$ is a location parameter, $\Sigma$ is a positive definite scatter matrix and $\nu$ is a degrees of freedom parameter. Here, $\nu$ is also known as shape parameter that controls the peakedness of density (see Jensen (1994)). When $\nu$ tends to infinity, the distribution will be $p$-variate normal distribution with mean vector $\boldsymbol{\mu}$ and covariance matrix $\Sigma$.

If $\boldsymbol{X} \sim t_p(\boldsymbol{\mu}, \Sigma, \nu)$ then the expectation and variance of $\boldsymbol{X}$ is given

$$E(\boldsymbol{X}) = \boldsymbol{\mu},$$
$$Var(\boldsymbol{X}) = \frac{\nu}{\nu - 2} \Sigma, \ if \ \nu > 2. \qquad (4)$$

The random vector $X$ from $p$-varaite t distribution has the following scale mixture representation

$$X = \mu + Y/\sqrt{U/\nu}, \tag{5}$$

where $Y$ is a $p$-variate normal random vector with mean $\mathbf{0}$ and covariance matrix $\Sigma$, $U$ is the chi-squared random variable with degrees of freedom $\nu$ and $U$ is independent of $Y$.

The conditional distribution of $X$ given $U = u$ is

$$X|u \sim N_p(\mu, \Sigma/u). \tag{6}$$

Using this conditional distribution, we obtain the following joint pdf of $X$ and $U$

$$f(x, u) = \frac{\left(\frac{\nu}{2}\right)^{\nu/2}}{\Gamma\left(\frac{\nu}{2}\right)(2\pi)^{p/2}} |\Sigma|^{-\frac{1}{2}} u^{\frac{\nu+p}{2}-1} e^{-\frac{u}{2}\left(\nu + (x-\mu)^T \Sigma^{-1}(x-\mu)\right)}. \tag{7}$$

Also, it can be easily shown that the conditional density function of $U$ given $X = x$ is

$$U|x \sim Gamma\left(\frac{\nu+p}{2}, \frac{1}{2}\left(\nu + (x-\mu)^T \Sigma^{-1}(x-\mu)\right)\right). \tag{8}$$

Further, using the conditional density function given in (8), we get the following conditional expectations will be used in the EM algorithm

$$E(U|x) = \frac{\nu + p}{\nu + (x-\mu)^T \Sigma^{-1}(x-\mu)}, \tag{9}$$

$$E(\log U \,|x) = \psi\left(\frac{\nu+p}{2}\right) - \log\left(\frac{1}{2}\left(\nu + (x-\mu)^T \Sigma^{-1}(x-\mu)\right)\right), \tag{10}$$

where $\psi(\alpha) = \frac{\Gamma'(\alpha)}{\Gamma(\alpha)}$ is the Digamma function.

### 3. ML estimation

Let $x_1, \ldots, x_n$ be a $p$-dimensional random sample from multivariate t-distribution with the parameters $\mu, \Sigma$ and $\nu$. The log-likelihood function of the multivariate t distribution is given by

$$\ell(\mu, \Sigma, \nu; x) = n \log \Gamma\left(\frac{\nu+p}{2}\right) - n \log \Gamma\left(\frac{\nu}{2}\right) - \frac{n\nu}{2} \log(\nu) - \frac{np}{2} \log \pi \\ - \frac{n}{2} \log|\Sigma|^{-1} - \frac{\nu+p}{2} \sum_{i=1}^{n} \log(\nu + s_i), \tag{11}$$

where $s_i = (x_i - \mu)^T \Sigma^{-1}(x_i - \mu)$. To find the ML estimators, the log-likelihood function is differentiated with respect to the parameters. These procedures yield the following ML estimating equations, we have to obtain the following score function for the parameters

$$\frac{\partial \ell(\mu, \Sigma, \nu; x)}{\partial \mu} = \sum_{i=1}^{n} (\nu + p) \frac{(x_i - \mu)}{\nu + s_i} = 0, \tag{12}$$

$$\frac{\partial \ell(\mu, \Sigma, \nu; x)}{\partial \Sigma^{-1}} = \frac{n}{2}\Sigma - \frac{\nu+p}{2}\sum_{i=1}^{n}\frac{(x_i - \mu)(x_i - \mu)^T}{\nu + s_i} = 0, \tag{13}$$

$$\frac{\partial \ell(\mu, \Sigma, \nu; x)}{\partial \nu} = \frac{1}{2}\sum_{i=1}^{n}\left(\psi\left(\frac{\nu+p}{2}\right) - \psi\left(\frac{\nu}{2}\right) - \log\nu - 1 - \log(\nu + s_i) - \frac{\nu+p}{\nu+s_i}\right) = 0. \tag{14}$$

Solving these equations yield the ML estimates for the parameters.

Concerning the local robustness properties of the ML estimators, the score functions should be bounded. However, when we examine the score functions for the t distribution given below,

$$\Psi_\mu = \frac{\nu+p}{\nu+s}\Sigma^{-1}(x-\mu), \tag{15}$$

$$\Psi_\Sigma = \frac{1}{2}\Sigma - \frac{\nu+p}{2}\frac{(x-\mu)(x-\mu)^T}{\nu+s}, \tag{16}$$

$$\Psi_\nu = \frac{1}{2}\log\nu + \frac{1}{2} + \frac{1}{2}\psi\left(\frac{\nu+p}{2}\right) - \frac{1}{2}\psi\left(\frac{\nu}{2}\right) - \frac{1}{2}\log(\nu+s) - \frac{1}{2}\frac{\nu+p}{\nu+s}, \tag{17}$$

we observe that the score function for $\mu$ tends to zero when $s = (x-\mu)^T\Sigma^{-1}(x-\mu)$ tends to $\infty$ and the score function for $\Sigma$ is bounded. However, the score function for $\nu$ tends to $-\infty$ when $s$ tends to $\infty$. The unboundedness of $\Psi_\nu$ can also be seen in Figure 1. The plot in Figure 1 is for the case $\mu = [0; 0]$, $\Sigma = I$ and $\nu = 3$.

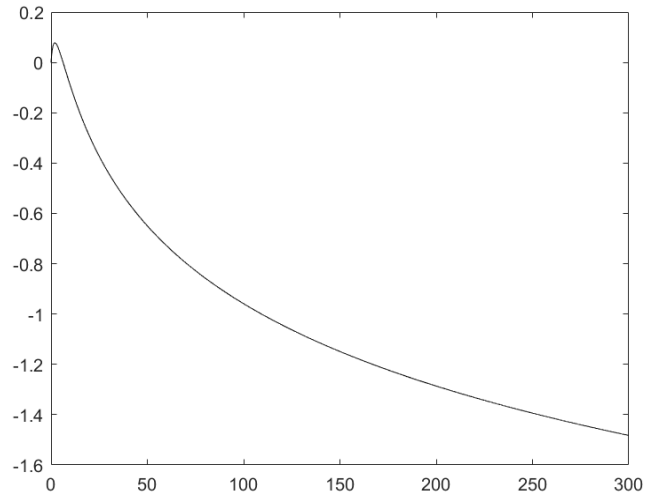

**Figure 1.** Score function plot of $\nu$ obtained from ML estimation.

Therefore, when we estimate the degrees of freedom parameter along with $\mu$ and $\Sigma$ the ML estimators will be no longer locally robust measured by the influence function. Note that the influence function for the ML estimators will be $J^{-1}\Psi$, $J$ is the Fisher information matrix and $\Psi$ is the score vector (see Hampel et al. (1986)). Thus, to get robust estimators in terms of influence function the degrees of freedom parameter is assumed to be known and other two parameters are estimated. When it is assumed that $\nu$ is fixed, the estimating equations for $\mu$ and $\Sigma$ will be as follows

$$\hat{\mu} = \frac{\sum_{i=1}^{n}\hat{w}_i x_i}{\sum_{i=1}^{n}\hat{w}_i}, \tag{18}$$

$$\hat{\Sigma} = \frac{1}{n} \sum_{i=1}^{n} \hat{w}_i \, (x_i - \hat{\mu})(x_i - \hat{\mu})^T, \qquad (19)$$

where $\hat{w}_i = (\nu + p) / (\nu + \hat{s}_i)$ and $\hat{s}_i = (x_i - \hat{\mu})^T \hat{\Sigma}^{-1}(x_i - \hat{\mu})$. Note that these equations can be viewed as adaptively sample mean and sample covariance matrix, where weights depend on the Mahalanobis distance between $x_i$ and $\mu$. The weight function $w(s) = (\nu + p) / (\nu + s)$ is a decreasing function of $s$ so that the outlying observations are downweighted by the corresponding weights (see Arslan et al. (1995), Kotz and Nadarajah (2004) and Nadarajah and Kotz (2008)).

To obtain the estimates the iteratively-reweighted algorithm which can be identified as an EM algorithm, can be used. In the following paragraph, we will describe the EM algorithm to obtain the ML estimators.

**EM algorithm to compute the ML estimates:**

Assuming that $U$ is missing in the scale mixture representation of the t distribution given in (5) we can implement the EM algorithm as follows. For the further details see McLachlan and Krishnan (1997). Using the conditional distribution given in (6) we have

$$\begin{aligned} X_i | u_i &\sim N_p(\mu, \Sigma/u_i), \\ U_i &\sim Gamma\left(\frac{\nu}{2}, \frac{\nu}{2}\right). \end{aligned} \qquad (20)$$

and using these equations, the complete data log-likelihood function can be obtained as

$$\begin{aligned} \ell_c(\Theta; x, u) = &\frac{n\nu}{2} \log\left(\frac{\nu}{2}\right) - n \log \Gamma\left(\frac{\nu}{2}\right) + \left(\frac{\nu + p}{2} - 1\right) \log u_i \\ &- \frac{np}{2} \log(2\pi) - \frac{n}{2} \log|\Sigma| - \frac{1}{2} \sum_{i=1}^{n} u_i \left(\nu + (x_i - \mu)^T \Sigma^{-1}(x_i - \mu)\right), \end{aligned} \qquad (21)$$

where $x = (x_1^T, \ldots, x_n^T)^T$ and $u = (u_1, \ldots, u_n)$. To handle the latency of $U$, we have to calculate the conditional expectation of the complete data log-likelihood function for given the observed data $x_i$

$$\begin{aligned} E(\ell_c(\Theta; x, u)|x_i) = &\frac{n\nu}{2} \log\left(\frac{\nu}{2}\right) - n \log \Gamma\left(\frac{\nu}{2}\right) + \left(\frac{\nu + p}{2} - 1\right) E(\log U_i | x_i) \\ &- \frac{np}{2} \log(2\pi) - \frac{n}{2} \log|\Sigma| - \frac{1}{2} \sum_{i=1}^{n} E(U_i | x_i)\left(\nu + (x_i - \mu)^T \Sigma^{-1}(x_i - \mu)\right). \end{aligned} \qquad (22)$$

The conditional expectations $E(U_i|x_i)$ and $E(\log U_i|x_i)$ can be calculated by using the conditional expectations given in (9) and (10). The steps of the EM algorithm can be given as follows:

**Steps of the EM algorithm:**

**1.** Set initial estimates $\Theta^{(0)} = \left(\mu^{(0)}, \Sigma^{(0)}, \nu^{(0)}\right)$ and a stopping rule $\epsilon$.
**2. E-Step:** Compute the following conditional expectations for $k = 0,1,2,\ldots$ iteration

$$\hat{u}_{1i}^{(k)} = E\left(U_i | x_i, \hat{\Theta}^{(k)}\right) = \frac{\hat{\nu}^{(k)} + p}{\hat{\nu}^{(k)} + (x_i - \hat{\mu}^{(k)})^T \hat{\Sigma}^{(k)^{-1}}(x_i - \hat{\mu}^{(k)})}, \qquad (23)$$

$$\begin{aligned} \hat{u}_{2i}^{(k)} &= E\left(\log U_i | x_i, \hat{\Theta}^{(k)}\right) \\ &= \psi\left(\frac{\hat{\nu}^{(k)} + p}{2}\right) - \log\left(\frac{1}{2}\left(\hat{\nu}^{(k)} + (x_i - \hat{\mu}^{(k)})^T \hat{\Sigma}^{(k)^{-1}}(x_i - \hat{\mu}^{(k)})\right)\right). \end{aligned} \qquad (24)$$

Then, we obtain the following objective function

$$Q(\Theta; \widehat{\Theta}^{(k)}) = \frac{n\nu}{2}\log\left(\frac{\nu}{2}\right) - n\log\Gamma\left(\frac{\nu}{2}\right) + \left(\frac{\nu+p}{2} - 1\right)\hat{u}_{2i}^{(k)}$$
$$- \frac{np}{2}\log(2\pi) - \frac{n}{2}\log|\Sigma| - \frac{1}{2}\sum_{i=1}^{n}\hat{u}_{1i}^{(k)}\left(\nu + (x-\mu)^T\Sigma^{-1}(x-\mu)\right). \quad (25)$$

**3. M-step 1:** Maximize $Q(\Theta; \widehat{\Theta}^{(k)})$ with respect to the unknown parameters $(\mu, \Sigma)$ to get the $(k+1)th$ estimates. This maximization gives the following equations

$$\hat{\mu}^{(k+1)} = \frac{\sum_{i=1}^{n}\hat{u}_{1i}^{(k)}x_i}{\sum_{i=1}^{n}\hat{u}_{1i}^{(k)}}, \quad (26)$$

$$\hat{\Sigma}^{(k+1)} = \frac{1}{n}\sum_{i=1}^{n}\hat{u}_{1i}^{(k)}\left(x_i - \hat{\mu}^{(k)}\right)\left(x_i - \hat{\mu}^{(k)}\right)^T. \quad (27)$$

**4. M-step 2:** If $\nu$ is estimated the following step should be implemented. Using the new values of $(\mu, \Sigma)$ which are obtained in M-Step 1, the following equation should be solved to get the $(k+1)th$ estimate for $\nu$

$$\sum_{i=1}^{n}\left(-\psi\left(\frac{\nu}{2}\right) + \log\left(\frac{\nu}{2}\right) + 1 + \hat{u}_{2i}^{(k)} - \hat{u}_{1i}^{(k)}\right) = 0. \quad (28)$$

**5.** Repeat E and M steps until the convergence rule $\|\widehat{\Theta}^{(k+1)} - \widehat{\Theta}^{(k)}\| < \epsilon$ is obtained.

## 4. ML*q* estimation

Let $x_1, \ldots, x_n$ be i.i.d random sample from multivariate t-distribution with the pdf given in (3). The ML$q$ estimators for the parameters of the multivariate t-distribution can be obtained by maximizing the following function

$$\ell_q = \sum_{i=1}^{n}L_q\big(f(x_i; \mu, \Sigma, \nu)\big), \quad q > 0, \quad (29)$$

where $0 < q < 1$ and $L_q$ function is given in (1). As $q \to 1$, we obtain the usual ML estimators. Taking the derivatives of $\sum_{i=1}^{n}L_q\big(f(x_i; \mu, \Sigma, \nu)\big)$ with respect to $(\mu, \Sigma, \nu)$, setting to zero and solving the resulting equations will give the ML$q$ estimators. These steps will give the following estimation equation

$$\sum_{i=1}^{n}U_q(x_i; \mu, \Sigma, \nu) = \sum_{i=1}^{n}U(x_i; \mu, \Sigma, \nu)f(x_i; \mu, \Sigma, \nu)^{1-q} = 0, \quad (30)$$

where $U(x_i; \mu, \Sigma, \nu) = \frac{\partial}{\partial\Theta}\log f(x_i; \mu, \Sigma, \nu)$ is the score vector and $\Theta = (\mu, \Sigma, \nu)$. After rearranging above equations for $\mu$ and $\Sigma$, we get

$$\hat{\mu}_q = \frac{\sum_{i=1}^{n}\widehat{w}_{qi}x_i}{\sum_{i=1}^{n}\widehat{w}_{qi}}, \quad (31)$$

$$\hat{\Sigma}_q = \frac{\sum_{i=1}^{n} \hat{w}_{qi}(x_i - \hat{\mu}_q)(x_i - \hat{\mu}_q)^T}{\sum_{i=1}^{n} \hat{v}_i}, \qquad (32)$$

where $\hat{w}_{qi} = \dfrac{\hat{v}_q + p}{(\hat{v}_q + \hat{s}_i)^{1 + \frac{(1-q)(\hat{v}_q + p)}{2}}}$, $\hat{v}_i = \dfrac{1}{(\hat{v}_q + \hat{s}_i)^{\frac{(1-q)(\hat{v}_q + p)}{2}}}$ and $\hat{s}_i = (x_i - \hat{\mu}_q)^T \hat{\Sigma}_q^{-1}(x_i - \hat{\mu}_q)$. Further, the ML$q$ estimator of $\nu$ can be found by solving the following equation

$$\sum_{i=1}^{n} \left(-\psi\left(\frac{\nu}{2}\right) + \psi\left(\frac{\nu + p}{2}\right) + \log \nu - \log(\nu + \hat{s}_i) - \hat{w}_i + 1\right) f(x_i; \hat{\mu}_q, \hat{\Sigma}_q, \nu)^{1-q} = 0. \qquad (33)$$

Note that $\hat{\mu}_q$ is similar to the $\hat{\mu}$ with slightly different weight function. For the $\hat{\Sigma}_q$ and $\hat{\Sigma}$ are different in terms of weighting. In the $\hat{\Sigma}_q$ instead of dividing $n$, we divide $\sum_{i=1}^{n} \hat{v}_i$, which makes more robust in terms of outliers. Concerning the parameter $\nu$, we observe that, unlike the ML case, the score function given in (33) is bounded as $s$ tends to $\infty$ provided that $q$ is finite and given. In Figure 2, which is for the case $\mu = [0; 0]$, $\Sigma = I$, $\nu = 3$ and $q = 0.85$, this behavior can be clearly noticed. Therefore, when we estimate the degrees of freedom along with $\mu$ and $\Sigma$ using ML$q$ estimation method, the resulting estimators will have bounded influence function which is not the case for ML estimators.

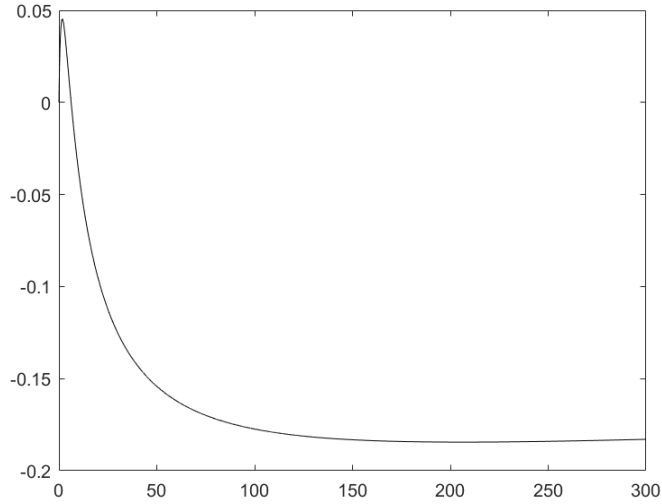

**Figure 2.** Score function plot of $\nu$ obtained from ML$q$ estimation.

Note that similar to the ML estimators, the ML$q$ estimators should be also computed using some numerical methods since the estimating equations cannot be solved analytically. To this extend, an EM-type algorithm similar to the algorithm proposed by Giuzio et al. (2016) for sparse and robust normal and t-portfolios by penalized L$q$-likelihood minimization will be proposed to obtain the ML$q$ estimates. The adaptation will be given in the following paragraph.

**EM-type algorithm to compute ML$q$ estimates:**

Let $\widehat{W}_{qi}^{(k)}(\hat{\mu}, \hat{\Sigma}, \hat{v}) = f(x_i; \hat{\mu}^{(k)}, \hat{\Sigma}^{(k)}, \hat{v}^{(k)})^{1-q}$. The estimates at $(k+1)th$ step will be obtained by maximizing the following function

$$\sum_{i=1}^{n} \widehat{W}_{qi}^{(k)}(\widehat{\boldsymbol{\mu}}, \widehat{\Sigma}, \hat{v}) \log f(\boldsymbol{x}_i, u_i; \boldsymbol{\mu}, \Sigma, v), \tag{34}$$

where $f(\boldsymbol{x}_i, u_i; \boldsymbol{\mu}, \Sigma, v)$ is the pdf of multivariate t-distribution. To implement the EM algorithm, we will use the complete data pdf $f(\boldsymbol{x}_i, u_i; \boldsymbol{\mu}, \Sigma, v)$ given in (7). For this case the estimating equation given in (30) can be rewritten as

$$\sum_{i=1}^{n} \widehat{W}_{qi}^{(k)}(\widehat{\boldsymbol{\mu}}, \widehat{\Sigma}, \hat{v}) U(\boldsymbol{x}_i, u_i; \boldsymbol{\mu}, \Sigma, v) = 0, \tag{35}$$

where $U(\boldsymbol{x}_i, u_i; \boldsymbol{\mu}, \Sigma, v) = \frac{\partial}{\partial \boldsymbol{\Theta}} \log f(\boldsymbol{x}_i, u_i; \boldsymbol{\mu}, \Sigma, v)$ is the derivative of the complete data log-likelihood function with respect to the parameters. Then, writing these derivatives in (35), we obtain the following equations for the parameters $\boldsymbol{\mu}, \Sigma$ and $v$

$$\sum_{i=1}^{n} u_i (\boldsymbol{x}_i - \boldsymbol{\mu}) W_{qi}(\widehat{\boldsymbol{\mu}}, \widehat{\Sigma}, \hat{v}) = 0, \tag{36}$$

$$\sum_{i=1}^{n} (\Sigma - u_i (\boldsymbol{x}_i - \boldsymbol{\mu})(\boldsymbol{x}_i - \boldsymbol{\mu})^T) W_{qi}(\widehat{\boldsymbol{\mu}}, \widehat{\Sigma}, \hat{v}) = 0, \tag{37}$$

$$\sum_{i=1}^{n} \left( \log\left(\frac{v}{2}\right) - \psi\left(\frac{v}{2}\right) + \log u_i - u_i + 1 \right) W_{qi}(\widehat{\boldsymbol{\mu}}, \widehat{\Sigma}, \hat{v}) = 0. \tag{38}$$

Then, to overcome the effect of latent variable on these equations we have to take conditional expectation of these equations for given $\boldsymbol{x}_i$. If we do so we get

$$\sum_{i=1}^{n} E(U_i | \boldsymbol{x}_i) (\boldsymbol{x}_i - \boldsymbol{\mu}) W_{qi}(\widehat{\boldsymbol{\mu}}, \widehat{\Sigma}, \hat{v}) = 0, \tag{39}$$

$$\sum_{i=1}^{n} (\Sigma - E(U_i | \boldsymbol{x}_i)(\boldsymbol{x}_i - \boldsymbol{\mu})(\boldsymbol{x}_i - \boldsymbol{\mu})^T) W_{qi}(\widehat{\boldsymbol{\mu}}, \widehat{\Sigma}, \hat{v}) = 0, \tag{40}$$

$$\sum_{i=1}^{n} \left( \log\left(\frac{v}{2}\right) - \psi\left(\frac{v}{2}\right) + E(\log U_i | \boldsymbol{x}_i) - E(U_i | \boldsymbol{x}_i) + 1 \right) W_{qi}(\widehat{\boldsymbol{\mu}}, \widehat{\Sigma}, \hat{v}) = 0. \tag{41}$$

We note that the conditional expectations $E(U_i|\boldsymbol{x}_i)$ and $E(\log U_i|\boldsymbol{x}_i)$ can be calculated by using the equations given in (9) and (10). Also, notice that we keep $W_{qi}(\boldsymbol{\mu}, \Sigma, v)$ this weight in the equations and called the resulting estimators as doubly reweighted estimators. Now, the steps of the EM algorithm can be as follows.

**EM-type algorithm for the MLq estimators:**

**1.** Take initial estimates $\boldsymbol{\Theta}^{(0)} = (\boldsymbol{\mu}^{(0)}, \Sigma^{(0)}, v^{(0)})$ and a stopping rule $\epsilon$.
**2. E-Step:** Calculate the conditional expectations for $k = 0,1,2,\ldots$ iteration

$$\hat{u}_{1qi}^{(k)} = E\left(U_i | \boldsymbol{x}_i, \widehat{\boldsymbol{\Theta}}^{(k)}\right) = \frac{\hat{v}_q^{(k)} + p}{\hat{v}_q^{(k)} + \left(\boldsymbol{x}_i - \widehat{\boldsymbol{\mu}}_q^{(k)}\right)^T \widehat{\Sigma}_q^{(k)^{-1}} \left(\boldsymbol{x}_i - \widehat{\boldsymbol{\mu}}_q^{(k)}\right)}, \tag{42}$$

$$\hat{u}_{2qi}^{(k)} = E\left(\log U_i | \boldsymbol{x}_i, \widehat{\boldsymbol{\Theta}}^{(k)}\right)$$

$$= \psi\left(\frac{\hat{v}_q^{(k)} + p}{2}\right) - \log\left(\frac{1}{2}\left(\hat{v}_q^{(k)} + \left(x_i - \hat{\mu}_q^{(k)}\right)^T \hat{\Sigma}_q^{(k)-1} \left(x_i - \hat{\mu}_q^{(k)}\right)\right)\right). \quad (43)$$

**3. M-step:** Write the conditional expectations given in E-step in solving equations (39)-(41) and rearrange these equations. Then, we obtain the following updating estimation equations

$$\hat{\mu}_q^{(k+1)} = \frac{\sum_{i=1}^n \hat{w}_{qi}^{(k)} x_i}{\sum_{i=1}^n \hat{w}_{qi}^{(k)}}, \quad (44)$$

$$\hat{\Sigma}_q^{(k+1)} = \frac{\sum_{i=1}^n \hat{w}_{qi}^{(k)} \left(x_i - \hat{\mu}_q^{(k)}\right)\left(x_i - \hat{\mu}_q^{(k)}\right)^T}{\sum_{i=1}^n \hat{v}_i^{(k)}} \quad (45)$$

where $\hat{w}_{qi}^{(k)} = \left(\hat{v}_q^{(k)} + p\right) \Big/ \left(\left(\hat{v}_q^{(k)} + \hat{s}_i^{(k)}\right)^{1 + \frac{(1-q)\left(\hat{v}_q^{(k)} + p\right)}{2}}\right)$, $\hat{v}_i^{(k)} = 1 \Big/ \left(\left(\hat{v}_q^{(k)} + \hat{s}_i^{(k)}\right)^{\frac{(1-q)\left(\hat{v}_q^{(k)} + p\right)}{2}}\right)$ and

$\hat{s}_i^{(k)} = \left(x_i - \hat{\mu}_q^{(k)}\right)^T \hat{\Sigma}_q^{(k)-1} \left(x_i - \hat{\mu}_q^{(k)}\right)$.

Use the following equation to obtain the new estimate for $v$

$$\sum_{i=1}^n \left(\log\left(\frac{v}{2}\right) - \psi\left(\frac{v}{2}\right) + \hat{u}_{2qi}^{(k)} - \hat{u}_{1qi}^{(k)} + 1\right) \widehat{W}_{qi}^{(k)}\left(\hat{\mu}_q^{(k)}, \hat{\Sigma}_q^{(k)}, v\right) = 0, \quad (46)$$

where $\widehat{W}_{qi}^{(k)}\left(\hat{\mu}_q^{(k)}, \hat{\Sigma}_q^{(k)}, v\right) = f\left(x_i; \hat{\mu}_q^{(k)}, \hat{\Sigma}_q^{(k)}, v\right)^{1-q}$.

**4.** Repeat these steps until convergence criteria $\left\|\widehat{\Theta}^{(k+1)} - \widehat{\Theta}^{(k)}\right\| < \epsilon$ is satisfied.

## 5. Simulation study

In this part, we will provide a simulation study to show the performance of the ML$q$ estimators over the ML estimators. The simulation study is performed by using MATLAB R2015b. The ML estimation procedure is done by using tdistfit in Matlab code for fitting multidimensional t-distributions (see the link https://github.com/robince/tdistfit for tdistfit code). For all numerical calculations, the stopping rule $\epsilon$ is taken as $10^{-6}$. Also, for the computations, we determine the following initial values of location, covariance matrix and degrees of freedom parameters for the ML and ML$q$ estimators. The initial values are taken as

$$\mu^{(0)} = mean(X), \ \Sigma^{(0)} = cov(X) \text{ and } v^{(0)} = 3,$$

where $X$ is the $p$-variate random sample from multivariate t distribution with location $\mu$, covariance matrix $\Sigma$ and degrees of freedom $v$.

We generate the data from multivariate t distribution using the stochastic representation given in (5) with the parameters

$$\mu = (\mu_{11}, \mu_{12})^T, \ \Sigma = \begin{bmatrix} \sigma_{1,11} & \sigma_{1,12} \\ \sigma_{1,21} & \sigma_{1,22} \end{bmatrix}, \ v.$$

We consider the following two cases

Case I : $\boldsymbol{\mu} = (2,1)^T$, $\Sigma = \begin{bmatrix} 1 & 0 \\ 0 & 1 \end{bmatrix}$, $\nu = 3$,

Case II : $\boldsymbol{\mu} = (2,1)^T$, $\Sigma = \begin{bmatrix} 2 & -0.5 \\ -0.5 & 2 \end{bmatrix}$, $\nu = 3$.

The tables include the mean and the mean Euclidean distance values of estimates, where the Euclidian distance of estimates are $\|\widehat{\boldsymbol{\mu}} - \boldsymbol{\mu}\|$ and $\|\widehat{\Sigma} - \Sigma\|$. We also note that the distance for $\hat{\nu}$ will be the mean squared error (MSE) which is given with the following formula

$$\widehat{MSE}(\hat{\nu}) = \frac{1}{N} \sum_{j=1}^{N} (\hat{\nu}_j - \nu)^2,$$

where $\nu$ is the true parameter value, $\hat{\nu}_j$ is the estimate of $\nu$ for the $jth$ simulated data. We take the sample sizes as $50, 100, 150$ and $200$. We set the replication number $(N)$ as $500$ for the simulation study. For the ML$q$ estimation, choosing $q$ is an important issue. In the simulation study, we choose $q$ which corresponds to minimum distance value of the mean Euclidian distance ($\|\widehat{\Theta} - \Theta\|$). In the simulation study, we investigate the behaviors of the estimators only for the outlier case. Without outliers, the estimators behave similar. We add extra five observations to the data. We consider adding five outliers as follows. Generate random numbers from uniform distribution and add them to the data.

Tables 1 and 2 display the simulation results for the sample sizes $50, 100, 150$ and $200$ with five outliers. We give mean and mean Euclidean distance values of estimates and true parameter values in tables for Cases I and II. We observe from the simulation results that the estimators for $\boldsymbol{\mu}$ have similar results. The estimators for $\Sigma$ obtained from ML$q$ seems slightly better than the estimators obtained from ML in terms mean Euclidian distance values. On the other hand, comparing the performance of estimators for the degrees of freedom parameter $\nu$, the ML$q$ estimator is definitely superior to the ML estimator in terms of MSE values. We observe that the estimates obtained from the ML$q$ are very close to the true values. This gets better when the sample sizes increases. For example, in Table 1 for the case $n = 200$, the mean of the estimated $\nu$s over the 500 replicates is 2.8369, which is very close to the true value $\nu = 3$ compare to the mean of the estimated $\nu$s obtained from ML method which is 1.7845. The similar results are noticed in Table 2 as well.

**Table 1.** Mean and mean Euclidean distance values of estimates for $n = 50, 100, 150$ and $200$ with the true parameter values given in Case I with five outliers.

|       |            |      | ML     |          | ML$q$  |          |
|-------|------------|------|--------|----------|--------|----------|
| $n$   | Parameter  | True | Mean   | Distance | Mean   | Distance |
| 50    | $\mu_{11}$ | 2    | 2.0181 | 0.2593   | 2.0119 | 0.2616   |
|       | $\mu_{12}$ | 1    | 1.0157 |          | 1.0099 |          |
|       | $\sigma_{1,11}$ | 1 | 1.1022 |          | 0.8674 |          |
|       | $\sigma_{1,12}$ | 0 | 0.2802 | 0.5670   | 0.0684 | 0.4177   |
|       | $\sigma_{1,22}$ | 1 | 1.1108 |          | 0.8800 |          |
|       | $\nu$      | 3    | 1.0861 | 3.6718   | 1.8827 | 1.3521   |
| 100   | $\mu_{11}$ | 2    | 2.0169 | 0.1773   | 2.0133 | 0.1794   |
|       | $\mu_{12}$ | 1    | 0.9989 |          | 0.9941 |          |
|       | $\sigma_{1,11}$ | 1 | 1.1120 |          | 0.9562 |          |
|       | $\sigma_{1,12}$ | 0 | 0.1476 | 0.3898   | 0.0327 | 0.3086   |
|       | $\sigma_{1,22}$ | 1 | 1.1274 |          | 0.9722 |          |
|       | $\nu$      | 3    | 1.4140 | 2.5294   | 2.4332 | 0.4458   |
| 150   | $\mu_{11}$ | 2    | 2.0095 | 0.1406   | 2.0064 | 0.1425   |
|       | $\mu_{12}$ | 1    | 1.0113 |          | 1.0085 |          |
|       | $\sigma_{1,11}$ | 1 | 1.1368 |          | 0.9944 |          |
|       | $\sigma_{1,12}$ | 0 | 0.1074 | 0.3379   | 0.0261 | 0.2581   |
|       | $\sigma_{1,22}$ | 1 | 1.1353 |          | 0.9945 |          |

| n | Parameter | True | Mean | Distance | Mean | Distance |
|---|---|---|---|---|---|---|
| | $\nu$ | 3 | 1.6255 | 1.9065 | 2.6853 | 0.1954 |
| 200 | $\mu_{11}$ | 2 | 2.0071 | 0.1226 | 2.0054 | 0.1235 |
| | $\mu_{12}$ | 1 | 0.9989 | | 0.9969 | |
| | $\sigma_{1,11}$ | 1 | 1.1498 | | 1.0177 | |
| | $\sigma_{1,12}$ | 0 | 0.0808 | 0.3054 | 0.0175 | 0.2238 |
| | $\sigma_{1,22}$ | 1 | 1.1384 | | 1.0045 | |
| | $\nu$ | 3 | 1.7845 | 1.4976 | 2.8369 | 0.1027 |

**Table 2.** Mean and mean Euclidean distance values of estimates for $n = 50, 100, 150$ and $200$ with the true parameter values given in Case II with five outliers.

| | | | ML | | ML$q$ | |
|---|---|---|---|---|---|---|
| $n$ | Parameter | True | Mean | Distance | Mean | Distance |
| 50 | $\mu_{11}$ | 2 | 2.0125 | 0.3659 | 2.0033 | 0.3741 |
| | $\mu_{12}$ | 1 | 1.0324 | | 1.0222 | |
| | $\sigma_{1,11}$ | 2 | 2.1293 | | 1.7309 | |
| | $\sigma_{1,12}$ | -0.5 | 0.0086 | 1.0750 | -0.2941 | 0.8972 |
| | $\sigma_{1,22}$ | 2 | 2.1153 | | 1.7226 | |
| | $\nu$ | 3 | 1.1135 | 3.5706 | 1.8607 | 1.3961 |
| 100 | $\mu_{11}$ | 2 | 2.0048 | 0.2470 | 2.0007 | 0.2508 |
| | $\mu_{12}$ | 1 | 1.0111 | | 1.0057 | |
| | $\sigma_{1,11}$ | 2 | 2.1851 | | 1.9206 | |
| | $\sigma_{1,12}$ | -0.5 | -0.2628 | 0.7110 | -0.4143 | 0.6205 |
| | $\sigma_{1,22}$ | 2 | 2.1891 | | 1.9296 | |
| | $\nu$ | 3 | 1.4705 | 2.3553 | 2.4789 | 0.3861 |
| 150 | $\mu_{11}$ | 2 | 2.0158 | 0.2046 | 2.0131 | 0.2065 |
| | $\mu_{12}$ | 1 | 0.9978 | | 0.9946 | |
| | $\sigma_{1,11}$ | 2 | 2.2093 | | 1.9615 | |
| | $\sigma_{1,12}$ | -0.5 | -0.3482 | 0.5883 | -0.4369 | 0.4959 |
| | $\sigma_{1,22}$ | 2 | 2.2169 | | 1.9686 | |
| | $\nu$ | 3 | 1.6632 | 1.8070 | 2.6878 | 0.1921 |
| 200 | $\mu_{11}$ | 2 | 1.9931 | 0.1757 | 1.9902 | 0.1778 |
| | $\mu_{12}$ | 1 | 1.0108 | | 1.0081 | |
| | $\sigma_{1,11}$ | 2 | 2.2916 | | 2.0335 | |
| | $\sigma_{1,12}$ | -0.5 | -0.4286 | 0.5800 | -0.4808 | 0.4555 |
| | $\sigma_{1,22}$ | 2 | 2.2879 | | 2.0284 | |
| | $\nu$ | 3 | 1.8045 | 1.4526 | 2.8067 | 0.1295 |

To further investigate the behavior of the estimators obtained from two methods we will simulate data from t distribution with the following parameter values

$$\boldsymbol{\mu} = (2,1)^T, \quad \Sigma = \begin{bmatrix} 1 & 0 \\ 0 & 1 \end{bmatrix}, \quad \nu = 2,$$

and some outliers generated from the uniform distribution. Then, the following table shows the estimated values obtained from two methods. Note that the difference between two sets of estimated values. Figure 3 displays the scatter plot of the data with the contour plots obtained from fitted densities. Note the fitted density obtained from the ML method is badly affected from the outlying observations at the right corner of the figure. Unlike the ML case, the ML$q$ fit seems resistant to these points.

**Table 3.** Estimation results for the first simulated data.

| Parameter | True | ML | MLq |
|---|---|---|---|
| $\mu_{11}$ | 2 | 2.1820 | 2.1139 |
| $\mu_{12}$ | 1 | 0.8844 | 0.9523 |
| $\sigma_{1,11}$ | 1 | 3.3389 | 1.7649 |
| $\sigma_{1,12}$ | 0 | -0.1126 | 0.1682 |
| $\sigma_{1,22}$ | 1 | 2.0851 | 1.4056 |
| $\nu$ | 2 | 1.4538 | 2.2446 |

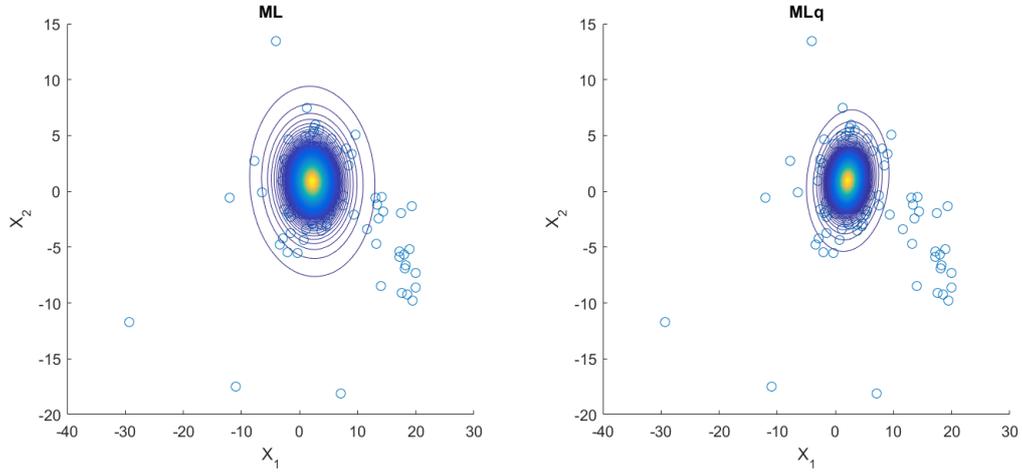

**Figure 3.** Scatter plot of the first simulated data along with the contour plots of the fitted densities obtained from ML and MLq.

Similar behavior is detected for the following example as well. In this case slightly different set of outliers are added to the simulated data.

**Table 4.** Estimation results for the second simulated data.

| Parameter | True | ML | MLq |
|---|---|---|---|
| $\mu_{11}$ | 2 | 2.2505 | 2.2488 |
| $\mu_{12}$ | 1 | 1.0159 | 1.0551 |
| $\sigma_{1,11}$ | 1 | 2.6206 | 1.2403 |
| $\sigma_{1,12}$ | 0 | -0.3807 | -0.1497 |
| $\sigma_{1,22}$ | 1 | 1.6122 | 1.0661 |
| $\nu$ | 2 | 1.3969 | 2.0636 |

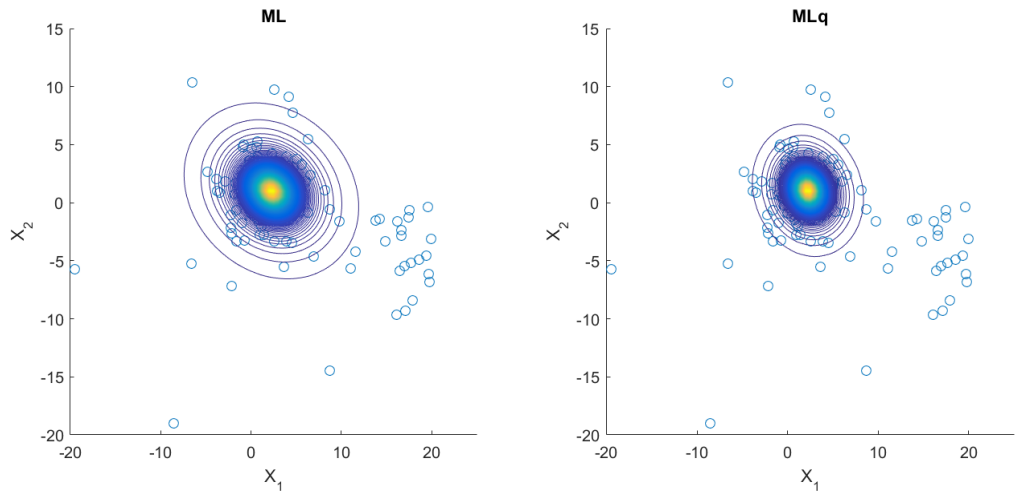

**Figure 4.** Scatter plot of the second simulated data along with the contour plots of the fitted densities obtained from ML and MLq.

## 6. Conclusions

In this paper, we have proposed the ML$q$ estimators for the parameters of the multivariate t distribution as an alternative to the ML estimators. We have provided an EM-type algorithm to compute the ML$q$ estimators. For the comparison, we have given a simulation study to illustrate performance of the proposed estimators over the ML estimators. We have observed from simulation results that the proposed method is working accurately to estimate all the parameters. Also, we see that the ML$q$ estimators outperform the ML estimators according to the mean Euclidian distance values for the parameters $\Sigma$ and $\nu$ and give similar results for the parameter $\boldsymbol{\mu}$ in the outlier case.